\documentclass{amsart}
 \usepackage{amsmath}
 \usepackage{amssymb}
 \usepackage{amscd}
 \usepackage[mathscr]{eucal}
\usepackage{epsfig}
\usepackage{graphics}

 \newtheorem{thm}{Theorem}[section]
 \newtheorem{lem}[thm]{Lemma}
 
 \newtheorem{example}[thm]{Example}
 \newtheorem{cor}[thm]{Corollary}
 \newtheorem{prop}[thm]{Proposition}
 \newtheorem{defn}[thm]{Definition}
 
 \newtheorem{conj}[thm]{Conjecture}
 \newtheorem{prob}[thm]{Problem}
 \newcommand{\B}{\mathscr{B}}

\newcommand{\Z}{\mathbb{Z}}

\begin{document}
 
 \title{Inequalities for the h- and flag h-vectors of geometric lattices}

\author{Kathryn Nyman \and Ed Swartz}
\maketitle

\begin{abstract}
  We prove that the order complex of a geometric lattice has a convex
  ear decomposition.  As a consequence, if $\Delta(L)$ is the order
  complex of a rank $(r+1)$ geometric lattice $L,$ then for all $i
  \le r/2$ the $h$-vector of $\Delta(L)$ satisfies $h_{i-1} \le h_i$
  and $h_i \le h_{r-i}.$
  
  We also obtain several inequalities for the flag $h$-vector of
  $\Delta(L)$ by analyzing the weak Bruhat order of the symmetric
  group.  As an application, we obtain a zonotopal {\bf cd}-analogue
  of the Dowling-Wilson  characterization of geometric lattices which
  minimize Whitney numbers of the second kind.  In addition, we are
  able to give a combinatorial flag $h$-vector proof of $h_{i-1} \le
  h_i$ when $i \le \frac{2}{7}(r+\frac{5}{2}).$
  
  \end{abstract}

\section{Introduction}
  The order complex of a geometric lattice is one of many simplicial
  complexes associated to matroids.  For a  geometric lattice $L,$ the
  order complex of $L, \Delta(L), $ is the simplicial complex whose
  simplices consist of all chains in $L, \hat{0} \neq x_1 < x_2 < .. <
  x_k \neq \hat{1}.$ The number of flats in each rank, also known as
  the Whitney numbers of the second kind, can be viewed as special
  cases of the flag $f$-vector  of $\Delta(L).$  The Euler
  characteristic of $\Delta(L)$ is the M\"{o}bius invariant of $L$
  \cite{Fo}.  Surveys of these topics are \cite{Ai}, \cite{Bj} and \cite{Za}. 
  
  Other enumerative invariants of $\Delta(L)$ have not received as
  much attention.  The explicit relationship between the flag
  $h$-vector of $\Delta(L)$ and the {\bf cd}-index of oriented
  matroids and zonotopes discovered in \cite{BER} suggests that it may
  be time to study the $h$-vector and flag $h$-vector of $\Delta(L).$  
  
  We begin with a review of the basic notions associated to geometric
  lattices, graded posets, $h$-vectors and flag $h$-vectors.  This is
  followed by an examination of the geometric lattices which minimize
  and maximize the flag $h$-vector of their order complex.  As a
  consequence we find the  zonotopes  with a specified dimension and
  number of zones which minimize or maximize the {\bf cd}-index.
  Then we show that $\Delta(L)$ has a convex ear decomposition.  An
  immediate consequence of this decomposition is our main theorem
  concerning the $h$-vectors of order complexes of geometric lattices.
  
  \begin{thm} \label{main h}
  Let $L$ be a rank $(r+1)$ geometric lattice.  Then, for $i \le r/2,$ the $h$-vector of $\Delta(L)$ satisfies

\begin{equation} \label{going up}
h_{i-1} \le h_i  \end{equation}
\begin{equation} h_i \le h_{r-i}
\end{equation}
\end{thm}

As is frequently the case with theorems of this type, the proof is not
combinatorial.  The remaining sections are devoted to understanding
the flag $h$-vector of $\Delta(L)$ with an eye toward providing a
combinatorial proof of (\ref{going up}).  This will lead us to an
examination of the weak Bruhat order on the symmetric group.
  
  \section{Definitions}

We will take all posets in this work to be finite.  A poset $P$ is
{\it graded} if all maximal chains have the same length and we call this
length the {\it rank} of $P$.  A graded poset has an associated
rank function $\rho$ which assigns to each element $y$ of $P$ a
positive integer such that $\rho (y) =k$, where $k$ is the length of
the longest chain of the form $y_0 < y_1 < \cdots < y_k = y$.

A {\it lattice} is a poset such that each pair of 
elements, $x$ and $y$, has a least upper bound, or {\it join}, 
denoted $x \vee y$, and a greatest lower
bound, or {\it meet}, denoted $x \wedge y$.  
Consequently, a lattice has a unique minimal element $\hat 0$ 
such that $x \geq \hat 0$ for
all $x \in L$, and a unique maximal element $ \hat 1$ with $x \leq 
\hat 1$ for all $x
\in L$.  We call an element of $L$ which covers $\hat 0$ an {\it atom},
and say that $L$ is {\it
atomic} if every element in $L$ can be written as the join of atoms.    

A {\it geometric lattice} is a graded atomic lattice whose rank 
function satisfies the
{\it semimodular condition} that for any $x$, $y \in L$,
\[\rho(x \vee y) + \rho(x \wedge y) \leq \rho(x) + \rho(y).\]

\begin{figure}[hbt]
\begin{center}
\input{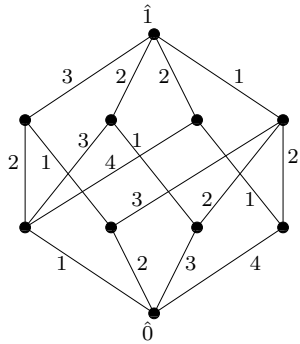}
\end{center}
\caption{An EL-labeling of the rank 3 near pencil on 4 atoms.}
\label{r-labeling}
\end{figure}

A broad class of geometric lattices arise from the affine dependencies
of a finite set of points $X$ in Euclidean space.  In this case the rank
$k$ elements of the lattice are subsets of the form $T \cap X$ where
$T$ is a $(k-1)$-dimensional subspace spanned by the elements of $X$.
These subsets are ordered by inclusion.  Points are in {\it general
position} if every set of  $k+1$ points spans a $k$-dimensional subspace.
One particularly useful geometric lattice arises from the {\it near pencil} 
arrangement on $n$ points in $r$-dimensional space which consists
of $(n-r+1)$ points on a line with the remaining $(r-1)$ points in
general position.  For ease of reference
we shall call this lattice the rank $r+1$ near pencil on $n$ 
atoms (see Figure \ref{r-labeling}).

\vspace{.5 cm}
 
A {\it matroid} $M = (X,\overline{\phantom{A}})$ is a set $X$ (for us
$X$ will always be finite) 
with a closure operation satisfying the exchange property (see 
\cite[Section 1.4]{Ox} for more details).  For $A \subseteq X$ denote the 
{\it closure of $A$} by $\overline{A}$.
A {\it simple matroid} is a matroid such that $\overline{\emptyset} =
\emptyset$ and $\overline{a} = a$ for every
element $a \in X$.
The closed sets, or {\it flats}, of a matroid, when partially ordered by 
inclusion form a
geometric lattice.  In fact a result of Birkhoff shows 
there is a bijection between geometric lattices and simple matroids \cite{Bi}.

A set $S\subseteq X$ is {\it independent} if $x
\notin \overline{S - x}$ for any $x \in S$.  
A {\it basis} of a matroid is a
maximal independent set and a {\it circuit} is a minimal dependent set.
A {\it loop} in a matroid $M=(X,\overline{\phantom{A}})$ is
an element $e \in X$ that is contained in no basis, while a {\it
coloop} is an element which is contained in every basis of the matroid.
Let $B$ be a basis of a matroid $M$.
If $e \notin B$, then $B \cup e$ contains a unique circuit,
$C(e,B)$.  The element $e$ is in $C(e,B)$, and we call
$C(e,B)$ the {\it fundamental circuit of e with respect to B}.  

\vspace{.5 cm}
By a basis of a geometric lattice we will mean a collection of atoms
whose cardinality is the rank of $P$ and whose join is $\hat 1$.  
Particular bases which will be
useful to us are the {\it nbc-bases}.  In order to define an
nbc-basis we first fix an arbitrary linear order $\omega$ on the
atoms of the geometric lattice $L$.  A {\it broken circuit} of
$(L,\omega)$ is a circuit with its least element removed.  The 
nbc-bases of $(L,\omega)$ are the bases of $L$ that do not contain a
broken circuit.  All such bases must contain the least atom.  Indeed,
if $B$ is a basis of $L$ which does not contain the least element,
then it will contain the broken circuit formed by removing the least
element from the fundamental circuit contained in the union of $B$ and
the least element.
\begin{example}
Let $L$ be the geometric lattice of the rank 3 matroid on
$\{1,2,3,4,5\}$ (with the natural order) whose bases consist of all
triples except $\{1,2,3\}$ and $\{3,4,5\}$.  Then the nbc-bases of
$L$ are the triples $\{1,2,4\}$, $\{1,2,5\}$, $\{1,3,4\}$,
$\{1,3,5\}$.  Notice that $\{1,4,5\}$ is not an nbc-basis because it
contains the broken circuit $\{4,5\}$.
\end{example}

\vspace{.3 cm}

Let $\Delta$ be a simplicial complex, i.e., 
$\Delta$ is a collection of subsets of a vertex set $X$
satisfying $x \in \Delta$ for any $x \in X$ and if $F \in \Delta$
and $G \subseteq F$ then $G \in \Delta$.
Maximal faces of $\Delta$ are
{\it facets} and $\Delta$ is {\it pure} if all its facets have
the same dimension.  A pure $d$-dimensional simplicial complex is said
to be {\it shellable} if there is an ordering of its facets $F_1, F_2,
\ldots, F_t$ such that $F_j \cap \bigcup_{i=1}^{j-1} F_i$ is a pure
$(d-1)$-dimensional complex for $j=2, \ldots, t$.  Such an ordering is
called a {\it shelling}.
Equivalently, a linear ordering $\psi$ on the facets of a complex is a
shelling if and only if it satisfies the following criterion 
(see, for instance,  \cite{Bj} )

\noindent {\bf Property $M$}:  For all facets $F$ and $F'$ of $\Delta$ such
that $F' <^{\psi} F$ there is a facet $F''$ with $F'' <^{\psi}F$ such
that $F' \cap F \subseteq F'' \cap F$ and $|F'' \cap F|= |F|-1$.

Given a poset $P$ with $\hat{0}$ and $\hat{1}$, the {\it order complex} $\Delta (P)$ of $P$ is the
simplicial complex whose vertices are the elements of 
$P-\{\hat{0},\hat{1}\}$ and whose
simplices are the chains of $P-\{\hat{0},\hat{1}\}$. Thus the facets correspond to maximal
chains.  We say that a poset $P$ is shellable if there
exists a shelling of $\Delta (P)$.
 
A poset $P$ admits an {\it R-labeling} 
if there is a map from the edges of $P$ to the positive integers (or
more generally to some partially ordered set)
such that in any interval
$[x,y] = \{z\in P: x \leq z \leq y\}$ of $P$ there is a unique 
saturated chain with increasing labels (known as a rising chain).  
For $y_2$ covering $y_1$ in $P$, denote the label on edge
$(y_1,y_2)$ by $\lambda(y_1,y_2)$.  Then a rising chain in $[x,y]$
is a maximal chain $x=y_0< y_1 < \cdots < y_k =y$ with
$\lambda(y_0,y_1) \leq \lambda(y_1,y_2) \leq \cdots \leq
\lambda(y_{k-1},y_k)$.  

An {\it EL-labeling} \cite{Bj2} of a poset is an R-labeling in 
which the unique
rising chain in an interval $[x,y]$ comes first lexicographically
among all of the chains in $[x,y]$.  An EL-labeling of a
geometric lattice can be obtained by labeling the
atoms $\{1, \ldots, n\}$ and labeling
the edge $(x,y)$ with the minimal atom $j$ such that $x \vee j = y$
(see Figure \ref{r-labeling}).  We will call this ordering the  {\it
minimal labeling} of the facets of $\Delta(L)$  and for a given
facet $F$ we denote its minimal label $\lambda(F).$ 
Notice with this labeling distinct chains have distinct labeling
sequences.  

A result of Bjorner shows that ordering the maximal chains of an 
EL-labeled poset lexicographically on the chain labels gives a
shelling of the associated order complex \cite[Theorem 2.3]{Bj2}.

\section{The flag $f$- and flag $h$-vectors of
geometric lattices}

For $P$ a rank $r+1$ poset, the number of simplices in $\Delta (P)$
of cardinality $k$ is 
denoted $f_k(\Delta(P))$, and $f(\Delta(P))= (f_0, \ldots, f_{r})$ is known as the
{\it $f$-vector} of the order complex.  The {\it $h$-vector} of
$\Delta(P)$  is defined as $h(\Delta(P)) = (h_0, \ldots, h_{r})$ where
\[\sum_{i=0}^r f_i(x-1)^{r-i} = \sum_{i=0}^r h_i x^{r-i}.\]

Given a rank $r+1$ poset $P$, and $S \subseteq [r]= \{1,2,\ldots, r\}$, let $P_S$ be the
rank selected subposet of $P$ defined by $P_S = \{x \in P : \rho(x)\in
S, x=\hat 0,\mbox{ or } x=\hat 1\}$.  The number of maximal
chains of $P_S$ is denoted $f_S(P)$; that is $f_S(P)$ counts the
number of chains in $P$ in which the ranks are the elements of $S$.  
The collection $\{f_S\}$,
$S\subseteq [r]$, is known as the {\it flag $f$-vector} of the
poset. The flag $f$-vector gives a natural refinement of the
$f$-vector of a poset's associated order complex as 
$f_i(\Delta(P)) = \sum_{|S|=i}f_S(P)$.

Let $\B_{r+1,n}$ denote
the rank $r+1$ truncated Boolean algebra on $n$ atoms.  $\B_{r+1,n}$ is
isomorphic to the rank $n$ Boolean algebra $\B_n$ for rank $i$, $i\leq
r$, and rank $r+1$ of $\B_{r+1,n}$ consists of
the maximal element $\hat 1$.   Every $(r+1)$-subset
of $[n]$ is a basis of $\B_{r+1,n}$ and every $i$-subset $i \leq r+1$ is
independent.  We shall see that the flag $f$-vector for rank $r+1$ 
geometric lattices on $n$ atoms is maximized by $\B_{r+1,n}$.
 
Dowling and Wilson \cite{DW1} proved that for rank $r+1$ geometric 
lattices with $n$ atoms the singleton flags $f_{\{i\}}, 1 \le i \le r,$ are
 minimized by the near pencil lattice.  
In her Ph.D. thesis \cite{Ny}, the first
author gave an explicit formula for $f_S$ of the near pencil lattice and
proved that this minimizes $f_S$, $S \subseteq [r]$ for all geometric 
lattices of rank $r+1$ with
$n$ atoms. In this section we prove the stronger result that $h_S$ is
minimized by the near pencil for all $S \subseteq [r]$.

\vspace{.5 cm}
Define the flag $h$-vector
$\{h_S(P)\}$, for $S \subseteq [r]$, by
\begin{equation}
h_S(P) =  \sum_{T \subseteq S} (-1)^{|S|-|T|} f_T (P).
\label{hdef}
\end{equation}
Again, the flag $h$-vector refines the $h$-vector since 
$h_i(\Delta(P)) = \sum_{|S|=i}h_S(P)$.  
An extremely useful combinatorial
interpretation of the flag $h$-vector is that it counts the number 
of chains with a specified descent set in an R-labeled poset.
We describe this interpretation below.

For a maximal chain $m:\hat 0=x_0<x_1< \cdots < x_{r+1} = \hat 1$, the
{\it descent set} of $m$ is $D(m)=\{i : \lambda(x_{i-1},x_i) >
\lambda(x_i, x_{i+1})\}$.
The following proposition can be found in \cite[Theorem 2.7]{Bj2} and
\cite[Theorem 3.13.2]{St3}.

\begin{prop}\cite{Bj2}\cite{St3}  \label{combinatorial h_S}
For $P$ a graded poset that admits an $R$-labeling, $h_S(P)$ is
the number of maximal chains of $P$ with labels having descent set
$S$.
\end{prop}

Since equation \eqref{hdef} can be inverted to give $ f_S(P) = 
\sum_{T\subseteq S} h_T (P)$ we see that $f_S(P)$ counts the number of
maximal chains in an $R$-labeled poset $P$ with descent set contained
in $S$.  

The positions of descents in maximal chains of a rank $r+1$ 
poset $P$ can be
encoded using the {\it {\bf ab}-index} $\Psi(P)$ of the poset 
which we describe presently.  
Assign to each chain in $P$ a word in the non-commuting variables ``{\bf a}'' and
``{\bf b}'' by assigning an {\bf a} if consecutive edge labels in the chain
increase and {\bf b} if the labels decrease.  Summing over all the chains
in $P$ gives the {\bf ab}-index
\[ \Psi(P) = \sum_{S \subseteq [r]} h_S(P) u_S \]
where the word $u_S = u_1 u_2 \cdots u_{r}$ is given by
$u_i = {\bf a}$ if $i \notin S$ and $u_i = {\bf b}$ if $i \in S$.

Through the course of this paper we will use the descent set of a
chain and its corresponding {\bf ab}-monomial interchangeably where we take
{\bf a} to mean an ascent and {\bf b} to indicate a descent in the chain
label.  $D(S)$ will indicate the set of all permutations with descent
set $S$ and $m(S)$ will refer to the {\bf ab}-monomial with descent set $S$.

\vspace{.5cm}
Next we consider the lattices which minimize and maximize the
flag $h$-vector.
\begin{prop}  Let $L$ be a rank $r+1$ geometric lattice with $n$
atoms.  Then for all $S \subseteq [r]$ we have $h_S(L) \leq
h_S(\B_{r+1,n})$.
\label{maxh}
\end{prop}
\begin{proof}

Place a linear ordering on the atoms of $L$ and $\B_{r+1,n}.$  We
construct an injection from the minimal labelings of $L$ to the minimal
labelings of $\B_{r+1,n}$ which preserves  descent sets.  This will
prove the proposition since $h_S$ is the number of minimal labelings
with descent set $S.$
                                                                                
Let $\lambda(F)= (\lambda_1,\dots,\lambda_{r+1})$ be a minimal label of
a facet $F$ in $\Delta(L).$  The uniqueness of rising chains in any
interval implies that $\lambda(F)$ is completely determined by the
initial segment $(\lambda_1,\dots,\lambda_i),$ where $\lambda_i$ is the
label preceding the last descent of $\lambda(F).$  For instance, if
$\lambda(F)=(3,6,9,2,7,12),$ then $(3,6,9),$ plus the knowledge that
$9$ precedes the last descent of the label, determines $\lambda(F).$
Since every subset of cardinality $i$ is a flat of $\B_{r+1,n}$ and
$\lambda_{i+1} < \lambda_i,$ there is a  unique facet $F^\prime$ of
$\Delta(\B_{r+1,n})$ such that $\lambda(F^\prime) =
(\lambda_1,\dots,\lambda_i,\dots)$ and whose final descent is
immediately after $\lambda_i.$ The function which takes $\lambda(F)$ to
$\lambda(F^\prime)$ is the required map.

\end{proof}

Since $ f_S(P) = 
\sum_{T\subseteq S} h_T (P)$ we immediately have the following result.
\begin{cor}
Among rank $r+1$ geometric lattices with $n$ atoms, $\B_{r+1,n}$
maximizes $f_S(P)$ for all $S\subseteq [r]$.
\end{cor}

We look again to the near pencil when considering the
lattice which minimize the flag $h$-vector.
\begin{lem}
Let $S \subseteq [r].$  Fix $i<r+1.$  The number of orderings
of $1,\dots,r+1$ such that $r+1$ comes before $i$ and has descent set $S$ is
independent of $i.$
\label{ordering}
\end{lem}
\begin{proof}
Suppose $r+1$ is the $k$th element of the permutation.  
There are ${r-1} \choose {r-k}$ ways to choose the elements (including
$i$) which appear after $r+1$.  Since $r+1$ effectively splits the
permutation into two smaller permutations we can consider the elements to
the right of $r+1$ as $1', 2', \ldots, (r+1-k)'$.  Let $D_r$ be the number of
ways to arrange $1',\ldots, (r+1-k)'$ consistent with descent set $S$ and
let $D_l$ be the number of ways to arrange the elements to the left of
$r+1$ consistent with descent set $S$.  Then
\[\sum_{k=1}^{r+1} {{r-1}\choose {r-k}} D_r \times D_l\]
is the number of permutations of $[r+1]$ such that $r+1$ appears before
$i$.  This number is independent of $i$.
\end{proof}

\begin{thm}
Let $L$ be a rank $r+1$ geometric lattice with $n$ atoms and let $P$
be the rank $r+1$ near pencil on $n$ atoms.  Then
$$h_S(L) \ge h_S(P).$$
\label{minh}
\end{thm}
\begin{proof}
Let $\{e_1,\dots, e_n\}$ and $\{f_1,\dots,f_n\}$ be the atoms of $P$
and $L$ respectively and let  $\Delta(P)$ and $\Delta(L)$ be the 
corresponding order complexes.
Order the atoms of $P$ so that $e_1,\dots,e_{r-1}$ are the coloops of
$P.$  What are the nbc-bases of $P?$   Any basis must contain the
coloops and $e_{r},$ the least atom in the non-trivial line.  On the
other hand, any basis of the form $\{e_1, \dots, e_{r}, e_i\}, i \ge
r+1$ is an nbc-basis of $P.$  Recall that $h_S(P)$ is the number of
orderings of the nbc-bases of $P$ for which the ordering is a minimal
labeling of the corresponding maximal chain of flats in $\Delta(P)$ and
the descent set is $S.$    For $i=r+1$ this is all orderings.  For $i>r+1$
an ordering is a minimal labeling if and only if $e_i$ comes before
$e_{r}.$
 
Now order the atoms of $L$ so that $f_1, \dots, f_{r+1}$ is a basis of $L.$
Clearly $B=\{f_1, \dots, f_{r+1}\}$ is an nbc-basis of $L,$ any ordering of
$B$ is a minimal labeling of the corresponding maximal chain, and the
contribution of these orderings to $h_S(L)$ is the same as the contribution
of all the orderings of $e_1,\dots,e_{r+1}$ to $h_S(P).$
 
For each $i >r+1$ we form a basis $B_i$ of $L$ as follows.  Let $C_i$ be the
fundamental circuit of $f_i$ with respect to $B.$  Let $B_i = B \cup
\{f_i\} - \{f_j\},$ where $f_j$ is the second highest element of $C_i.$
For instance, if $B = \{f_1,f_2,f_3,f_4,f_5\}$ and the fundamental circuit
of $f_8$ is $C_8 = \{f_1, f_3,f_4,f_8\},$ then $B_8 =
\{f_1,f_2,f_3,f_5,f_8\}.$ Now one can check that each $B_i$ is an
nbc-basis of $L$ and that any ordering of $B_i$ where $f_i$ comes before
the least element of $C_i$ is a minimal labeling of the corresponding
maximal chain.  In combination with Lemma \ref{ordering} this
shows that there
are at least as many minimal labelings with descent set $S$ for
$\Delta(L)$ as there are for $\Delta(P).$
\end{proof}
When the number of atoms of a rank $r+1$ geometric lattice is not
specified the standard Boolean algebra $\B_{r+1}$ minimizes the flag
$h$-vector \cite[Proposition 7.4]{BER}.
\begin{prop}\cite{BER}  Let $L$ be a geometric lattice of rank $r+1$.
Then for all $S \subseteq [r]$ we have $h_S(L) \geq h_S(\B_{r+1})$.
Hence the {\bf ab}-index $\Psi(L)$ is coefficient-wise greater than or
equal to the {\bf ab}-index of the Boolean algebra $\B_{r+1}$.
\label{7.4ber}
\end{prop}

\vspace{.5cm}

Oriented matroids are signed versions of standard matroids.  We refer
the reader to \cite{Bj4} for more details.  The elements of the
oriented matroid, when partially ordered, form an Eulerian poset (see
\cite{St3}) which is known as the lattice of regions.
A poset is {\it Eulerian} if every interval $[x,y]$, where $x \not= y$, has
the same number of elements of odd rank as even rank.  
Of interest to our work is the fact that underlying each oriented
matroid is a standard matroid along with its associated geometric lattice of flats.   

A collection of hyperplanes $\mathscr{H} = \{H_e\}_{e\in E}$ is {\it
essential} if $\bigcap_{e\in E} H_e = \{0\}$.
A special class of oriented matroids, called realizable matroids, have an
associated essential hyperplane arrangement.  The lattice of regions of the
realizable oriented matroid is isomorphic to the face lattice of the
corresponding hyperplane arrangement.  Every essential hyperplane arrangement has
an associated {\it zonotope}, which is the polytope formed by taking
the Minkowski sum of the 
normals to the hyperplanes (see \cite{Bj4}).  

It was noted by Fine and proved by Bayer and Klapper \cite{BK} that
when $P$ is an Eulerian poset the {\bf ab}-index of $P$, $\Psi(P)$, can be
written in terms of ${\bf c} = {\bf a+b}$ and ${\bf d}={\bf a} \cdot
{\bf b }+ {\bf b} \cdot {\bf a}$.  
When $\Psi(P)$ is
expressed in terms of {\bf c} and {\bf d} it is referred to as the {{\bf cd}-index}.  
The face lattice of an $r$-dimensional convex polytope is an Eulerian
poset and the ${\bf cd}$-index of the polytope is defined to be the
{\bf cd}-index of its corresponding face lattice.

When $P$ is the lattice of regions of an oriented
matroid every occurrence of {\bf d} in the {\bf cd}-index appears as {\bf 2d} and
so it is referred to as the {\bf c}-{\bf 2d}-index \cite{BER}.
The lattice of flats of the underlying matroid
contains all of the information necessary to determine the {\bf c}-{\bf 2d}-index of
the lattice of regions of an oriented matroid \cite[Theorem 3.4]{BS}. 
This connection is made explicit in \cite{BER}.  The following
proposition indicates how to construct the {\bf c}-{\bf 2d}-index of a zonotope
given the geometric lattice underlying the associated hyperplane arrangement
\cite[Corollary 3.2]{BER}.

\begin{prop}\cite{BER}
Let $L$ and $Z$ be the underlying geometric lattice and zonotope
associated to an essential hyperplane arrangement respectively.  Then the
{\bf c}-{\bf 2d}-index of the zonotope $Z$ is given by 
\[\Psi(Z)=\omega({\bf a} \cdot
\Psi(L))\] 
where $\omega$ is a linear function which takes {\bf ab} words to {\bf cd}
words by replacing each occurrence of {\bf ab} with {\bf 2d}, then replacing
the remaining letters with {\bf c}'s.
\label{3.1ber}
\end{prop}

The hyperplane arrangement $\{x \in \mathbb{R}^{r+1} : x_i =0\}$ for
$i = \{1,\ldots, r+1\}$ has the Boolean algebra $B_{r+1}$ as its 
underlying lattice of flats.  The zonotope corresponding to this 
arrangement is the $(r+1)$-dimensional cube.  Propositions 
\ref{7.4ber} and \ref{3.1ber}
together imply the following \cite[Corollary 7.6]{BER}.

\begin{cor}\cite{BER}
Among all zonotopes of dimension $r+1$, the $(r+1)$-dimensional cube has
the smallest {\bf c}-{\bf 2d}-index.
\end{cor}

Consider the rank $r+1$ near pencil on $n$ atoms and the truncated
Boolean algebra, $\B_{r+1,n}$.  As seen in Section 2 the near pencil
is associated with the arrangement of $(n-r+1)$ points on a line
with the remaining $r-1$ points in general position.  Similarly
$\B_{r+1,n}$ can be associated to the arrangement of $n$ points in
general position in $\mathbb{R}^{r+1}$.  In both of these point
arrangements we can consider the set of rays from the origin to each
point.  Taking the hyperplanes normal to these rays gives an essential
hyperplane arrangement whose underlying geometric lattice is the near
pencil, and $\B_{r+1,n}$ respectively.
Combining Proposition \ref{3.1ber} with Theorem \ref{minh} and Proposition
\ref{maxh} gives the following analogous result for $(r+1)$-dimensional
zonotopes with $n$ zones.  

\begin{cor}  Let $\mathscr{H}_{P}$ and $\mathscr{H}_{\B_{r+1,n}}$ denote
essential hyperplane arrangements whose underlying geometric lattice
is the rank $r+1$ near pencil and truncated Boolean algebra on $n$
atoms respectively.
Among all zonotopes of dimension $r+1$ with $n$ zones, the zonotope
corresponding to $\mathscr{H}_{\B_{r+1,n}}$ has the
largest {\bf c}-{\bf 2d}-index, and the zonotope corresponding to $\mathscr{H}_P$
has the smallest {\bf c}-{\bf 2d}-index.
\end{cor}

\section{Convex ear decompositions}

Convex ear decompositions were introduced by Chari
\cite{Ch}.  Our convex ear decomposition is motivated by a basis for $H_\star(\Delta(L),\Z)$ constructed by Bj{\"o}rner \cite{Bj5}.
A {\it convex ear decomposition} of a pure
$(r-1)$-dimensional simplicial complex $\Delta$ is an
ordered sequence $\Delta_1,\Delta_2,\dots,\Delta_m$ of
pure $(r-1)$-dimensional subcomplexes of $\Delta$ such
that 

\begin{enumerate}
 \item
  $\Delta_1$ is the boundary complex of a simplicial
$r$-polytope, while for each $j=2,\dots,m,
\Delta_j$ is an $(r-1)$-ball which is a proper subcomplex
of the boundary of a simplicial $r$-polytope.

  \item
    For $j \ge 2, \Delta_j \cap (\bigcup ^{j-1}_{k=1}
\Delta_k) = \partial \Delta_j.$
 \item
$\bigcup^{m}_{k=1} \Delta_k = \Delta$.

\end{enumerate}

\begin{thm} \cite{Ch} \label{convex h}
  If $\Delta$ has a convex ear decomposition, then for
$i \le  r/2 $ the
$h$-vector of $\Delta$ satisfies

$$
 h_{i-1} \le h_i$$
 $$  h_i \le h_{r-i}
$$\end{thm}

For the rest of this section $(b_1,\dots,b_{r+1})$ is
an
ordered basis of atoms in $L$ with corresponding chain
of flats
$\hat{0}< x_1< \dots < x_r < \hat{1}, x_i = b_1 \vee \dots \vee b_i.$

\begin{lem}
  Let $(b_1,\dots,b_{r+1})$ be a minimal labeling of
a
facet of $\Delta(L).$  Then $B=\{ b_1,\dots,b_{r+1}\}$ is
an
nbc-basis of $L.$
\end{lem}

\begin{proof}
  Suppose $B$ contains
a broken circuit.   If $x_i$ is the lowest ranked flat which contains
a broken circuit, then $b_i \neq
\lambda(x_{i-1},x_i).$
\end{proof}

\begin{lem}[Switching lemma]
\label{switching lemma}
  Let $(b_1,\dots, b_i, b_{i+1},
\dots,b_{r+1})$ be a
minimal labeling of a facet of $\Delta(L).$ If
$b_i <
b_{i+1},$ then  $(b_1,\dots,  b_{i+1}, b_i,\dots,b_{r+1})$
is
also a minimal labeling of a facet of $\Delta(L).$
\end{lem}

\begin{proof}
  For two flats $y < x$ in $L$ let $\{x-y\} = \{\mbox{atoms }e: y \vee e = x\}.$ Suppose $(b_1,\dots,  b_{i+1}, b_i,\dots,b_{r+1})$ is
not a minimal labeling.  Then there exists an atom $e$
such that either $e \in \{(x_{i-1} \vee b_{i+1}) - x_{i-1}\}$ and $e < b_{i+1},$ or $e \in \{x_{i+1} -(x_{i-1} \vee b_{i+1})\} $ and $e <
b_i.$  In the first case, $e \in \{x_{i+1} - x_i\}.$ 
However, this implies that $b_{i+1}$ is not the least
atom in $x_{i+1} - x_i.$  In the second case, either
$e \in x_i - x_{i-1},$ which implies that $b_i$ is not
the minimal atom in $x_i - x_{i-1},$ or $e \in x_{i+1}
- x_i.$ But this last is impossible since $e < b_i <
b_{i+1}$ and $b_{i+1}$ is the least atom in $x_{i+1} -
x_i.$ 
 \end{proof}
  
Let $B$ be a basis of $L.$ Associated to any ordering $(b_1,\dots,b_{r+1})$ of $B$ is the facet $F =b_1< b_1 \vee b_2 < \dots < b_1 \vee  \dots \vee b_r$ of $\Delta(L).$ The {\it basis labeling} of $F$ (with respect to $B$) is $(b_1,\dots,b_{r+1}).$  This may or may not be the same as $\lambda(F).$  

Let $B_1,\dots,B_m$ be the nbc-bases of $L$ in lexicographic order.
For each $j, 1 \le j \le m,$ let $\Sigma_j$ be the union of all the
facets of $\Delta(L)$ associated to all possible orderings of $B_j.$
Each $\Sigma_j$ is isomorphic to the order complex of the rank $r+1$
Boolean algebra and as a simplicial complex is the boundary of the
first barycentric subdivision of the $r$-simplex.  Now define
$\Delta_j$ to be the pure subcomplex of $\Sigma_j$ whose facets are
the facets of $\Sigma_j$ whose minimal labeling and basis labeling
coincide.  Except for $\Sigma_1 = \Delta_1,$ each $\Delta_j$ is a
proper subcomplex of $\Sigma_j$ \cite[Lemma 7.6.2]{Bj}.  

\begin{prop}
  If $2 \le j \le m,$ then $\Delta_j$ is a closed $(r-1)$-ball.  
\end{prop}

\begin{proof}
It is sufficient to show that $\Delta_j$ is nonempty and shellable.  To
see that $\Delta_j$ is nonempty, we note that for any basis $B$ of $L$
the minimal ordering and the basis ordering are the same for the
maximal chain corresponding to ordering $B$ in reverse if and only if
$B$ is an nbc-basis.
                                                                                
Order the facets of $\Delta_j$  in {\it reverse} lexicographic order
with respect to the basis labeling of the corresponding maximal chains.
   We show that this ordering satisfies property $M.$  Suppose
$F^\prime$ and $F$ are  facets of  $\Delta_j$ and  $\lambda(F^\prime) <
\lambda(F).$ Let $c=(\hat{0}=x_0< x_1 < \dots < x_k<x_{k+1}=\hat{1})$
be the chain which represents their intersection.  We must find
$F^{\prime\prime},$ a facet of $\Delta_j,$ lexicographically after $F$
such that $c \subseteq (F \cap F^{\prime\prime})$ and $|F \cap
F^{\prime\prime}| = r-1.$ Let $m$ be the least index such that the
length of  $[x_m,x_{m+1}]$ is greater than one.  Since $\lambda(F)$ is
lexicographically before $\lambda(F^\prime),$ there are  basis atoms $b
< \hat{b}$ in $B_j$ such that the chain corresponding to $F$ contains
as a short saturated chain $y < (y \vee b) < (y \vee b \vee \hat{b})$
with $x_m \le y$ and $(y \vee b \vee \hat{b}) \le x_{m+1}.$  Let
$F^{\prime\prime}$ be the facet corresponding to interchanging $b$ and
$\hat{b}$ in $\lambda(F).$ Then $|F \cap F^{\prime\prime}| = r-1$ and
$\lambda(F^{\prime\prime}) < \lambda(F).$  By the switching lemma,
$F^{\prime\prime}$ is a facet of $\Delta_j.$
 \end{proof}

\begin{prop}
  If $j \ge 2,$ then $ \Delta_j \cap (\bigcup ^{j-1}_{k=1}
\Delta_k) = \partial \Delta_j.$
\end{prop}

\begin{proof}
   Let $G$ be a face in $ \Delta_j \cap (\bigcup ^{j-1}_{k=1}
\Delta_k)$. By definition $G$ is not a facet.  The boundary of $\Delta_j$ is equal to the boundary of $\overline{\Sigma_j - \Delta_j}$ (topological closure). So it is sufficient to show that $G$ is contained in a facet of $\overline{\Sigma_j - \Delta_j}$. 
   
Write $G = x_1 < \dots < x_k.$  By assumption $ G \subset F, F$ a
facet of $\Delta_j,$ and $G \subset F^\prime, F^\prime$ a facet whose
corresponding minimal labeling basis $B^\prime$ is lexicographically
before $B_j.$  Therefore, there is some pair $x_m, x_{m+1}$ such that
$B^\prime \cap \{x_{m+1} - x_m\}$ is lexicographically before $B_j
\cap \{x_{m+1} - x_m\}.$  Hence, the unique increasing minimally
labeled saturated chain of $[x_m, x_{m+1}]$ is not contained in $B_j.$
Now let $\hat{F}$ be a facet of $\Sigma_j$ obtained as follows.  First
saturate the interval $[x_m, x_{m+1}]$ by adding in the atoms of $B_j
\cap \{x_{m+1} - x_m\}$ in increasing order.  Then extend this to a
saturation of the chain corresponding to $G$ in any way which results
in  a facet of $\Sigma_j.$ Such an $F^{\prime\prime}$ contains $G$ and
must be in   $\overline{\Sigma_j - \Delta_j}$ since its minimal label
and its basis label are not equal.
\end{proof}

The two previous propositions show that $\Delta(L)$ has a convex ear
decomposition.  An immediate consequence is Theorem \ref{main h}.  
Since $h_0 \le h_1 \le \dots \le h_{\lceil r/2 \rceil},$ a natural
question is whether or not the $g$-vector of $\Delta(L)$ is an
$M$-vector.  The $g$-vector is $(g_0,g_1,\dots, g_{\lceil r/2
  \rceil}),$ where $g_i = h_i - h_{i-1}.$  A sequence of nonnegative
integers is an $M$-vector if it is the Hilbert function of a quotient
of a polynomial ring.  See, for instance, \cite[Theorem 2.2,
pg. 56]{St} for an equivalent numerical definition of $M$-vector.  

\begin{prob}
  Is the $g$-vector of the order complex of a geometric lattice an $M$-vector?
\end{prob}

Note: After this
paper was written the second author discovered a proof
that the g-vector of any space with a convex ear
decomposition is an M-vector.

\section{The weak Bruhat order}
  Let $(b_1,\dots,b_{r+1})$ be an nbc-basis of $L$ ordered so that
  $b_1 < \dots < b_{r+1}.$  Then we can identify all the orderings of
  the basis with $S_{r+1},$ the symmetric group on $r+1$ letters. For
  $\pi \in S_{r+1}$ we write $\pi = a_1 a_2 \cdots a_{r+1},$ where
  $a_i = \pi(i).$   The switching lemma tells us that if $\pi $
  corresponds to a minimal labeling and $\pi^\prime$ is obtained from
  $\pi$ by interchanging $a_i$ with $a_{i+1}$ when $a_i < a_{i+1},$
  then $\pi^\prime$ also corresponds to a minimal labeling.
  
\begin{defn}
 Let $\pi, \pi^\prime \in S_{r+1}.$  Then $\pi \le_w \pi^\prime$ if
 and only if $\pi^\prime$ can be obtained from $\pi$ by repeated
 application of the above switching procedure.
\end{defn}

Evidently $\le_w$ is a partial order.  It is, in fact, the  weak
Bruhat order on $S_{r+1}.$   An equivalent definition of $\le_w$  is
the following.  Define the {\it inversion} set of $\pi$ to be $I(\pi)
= \{ (a_i,a_j): a_i > a_j$ and $i < j\}.$  Then $\pi \le_w \pi^\prime$
if and only if $I(\pi) \subseteq I(\pi^\prime$).  See \cite{Bj3} for
more information on the weak order in general Coxeter groups.  

 Let $T,S \subseteq [r].$  We say $S$ {\it dominates} $T$ if there
 exists an injection $\phi:D(T) \to D(S)$  such that $\pi \le_w
 \phi(\pi)$ for all permutations $\pi \in D(T).$
 
  \begin{prop}
 If $S$ dominates $T,$ then $h_T \le h_S$ for all geometric lattices of rank $r+1$ or greater.
 \end{prop}
 \begin{proof}
 The switching lemma and Proposition \ref{combinatorial h_S} imply
 that there are at least as many facets which contribute to $h_S$ as
 $h_T$ at each step in the convex ear decomposition.
 \end{proof}

 How do we find pairs $T,S$ such that $S$ dominates $T?$  First some elementary facts.  
 \begin{prop} \label{lengthen}
    Suppose $S$ dominates $T.$ Let $u,v$ be   {\bf ab}-monomials,
    possibly equal to  $\emptyset$. Then
   \begin{enumerate}
    \item
      $m(S) \cdot {\bf a} \cdot v$ dominates $m(T) \cdot {\bf a} \cdot v,$

    \item
      $u \cdot {\bf a} \cdot m(S)$ dominates $u \cdot {\bf a} \cdot m(T),$
      
      \item
        $u \cdot {\bf a} \cdot m(S) \cdot {\bf a} \cdot v$ dominates
        $u \cdot {\bf a} \cdot m(T) \cdot  {\bf a} \cdot v.$
\end{enumerate}
\end{prop}

\begin{proof}
  We prove (3) since the other proofs are virtually identical.  By
  definition, if $\pi \le_w \pi^\prime \in S_{r+1},$ then $\pi(1) \le
  \pi^\prime(1)$ and $\pi(r) \ge \pi^\prime(r).$  Let $\phi:D(T) \to
  D(S)$ be an injection which preserves the weak Bruhat  order.  Let
  $$\pi=s_1 \cdots s_m a_1 a_2 \cdots a_{r+2} a_{r+3} t_1 \cdots t_k
  $$ be a permutation such that the {\bf ab}-monomial of  the $s_1
  \cdots s_m a_1$ is $u,$ the {\bf ab}-monomial of $a_1 \cdots
  a_{r+3}$ is ${\bf a} \cdot m(T) \cdot {\bf a},$ and the {\bf
  ab}-monomial of $a_{r+3 } \cdot t_1 \cdots t_k$ is $v.$   Identify
  the  ordered set $[r+1]$ with the  ordered set $[r+3] -
  \{a_1,a_{r+3}\}$ in the canonical way.  Define $\psi(\pi)$ to be the
  permutation  obtained by applying  $\phi$ to $a_2 \cdots a_{r+2}$
  using this identification.  Clearly, $\psi$ is an injection and
  $\pi \le_w \psi(\pi).$ Since $a_1 < a_2$ and $a_{r+2} < a_{r+3},$
  the descent monomial of $\psi(\pi)$ is $u \cdot {\bf a} \cdot m(S)
  \cdot {\bf a} \cdot v.$
\end{proof}

\begin{prop} \label{subset condition}
  If $S$ dominates $T,$ then $T \subseteq S.$  
\end{prop}

\begin{proof}
  Suppose $i \in T.$  Let $\pi$ be a permutation with descent set $T$
  such that $\{\pi(1),\dots,\pi(i)\} = \{r-i+2,\dots,r+1\}.$ If $\pi
  \le_w \pi^\prime,$ then $\{\pi^\prime(1),\dots,\pi^\prime(i)\}$ must
  also equal $\{r-i+2,\dots,r+1\}.$ Hence $\pi^\prime$ also has a
  descent at $i.$ 
  \end{proof}
  
Another place to look for $T \subseteq S$ with $S$ dominating $T$ is
through the symmetries of the flag $h$-vector of the Boolean algebra.
Since $h_T \le h_S$ for all geometric lattices of rank at least $r+1,$
we can begin our search by examining $\B_{r+1}.$  This lattice has
$\Z_2 \oplus \Z_2$ symmetry. Let $\beta \in S_{r+1}$ be the
permutation which reverses order, i.e., $\beta(i) = r-i+2.$  We omit
the elementary proof of the following:

\begin{prop}
  Let $T$ be the descent set of $\pi \in S_{r+1}.$
  \begin{enumerate}
   \item
     The descent set of $\beta \circ \pi$ is $[r]-T.$
     
    \item
     The descent set of $\pi \circ \beta$ is $T \circ \beta,$ where $T
     \circ \beta = \{i \in [r]: r-i+1 \notin T\}.$
    \item
     The descent set of $\beta \circ \pi \circ \beta$ is $r - [T],$
     where $r-[T] = \{i:r-i+1 \in T\}.$
 \end{enumerate}
\end{prop}

\noindent  Combined with Proposition \ref{combinatorial h_S}, the
     above proposition shows that in $\B_{r+1}, h_T = h_{r-T} =
     h_{[r]-T} = h_{T \circ \beta}.$  Proposition \ref{subset
     condition} rules out the possibility of $T$ being dominated by
     $r-T$ or $[r]-T$ except when they are equal.  However, if for
     each $i$ at most one of $i$ and $r-i+1$ is in  $T,$ then $T
     \subseteq T \circ \beta.$

\begin{example}
  Let $T= \{1\} \subseteq [3].$  Then $T \subset T \circ \beta =
  \{1,2\}.$  As the map which sends $\{2134\} \to \{3214\}, \{3124 \}
  \to \{4312\},$ and $ \{4123\} \to \{4213\}$ shows, $T$ is
  dominated by $T \circ \beta$ and hence $h_T \le h_{T \circ \beta}$
  for all geometric lattices of rank 4 (or more). 
\end{example}

\begin{conj} \label{main conjecture}
  If $T \subseteq T \circ \beta,$ then $T \circ \beta$ dominates $T.$
\end{conj}

We have verified this conjecture by computer for $r \le 8.$ For
convenience, Table \ref{known cases} lists all the cases that we have
used in Section \ref{i to i+1}. These computations were helped by the
observation that if $T \subseteq T \circ \beta$ and $T \circ \beta$
dominates $T,$ then $r-(T \circ \beta)$ dominates $r-T.$  To see this,
use the fact  that if $\pi \le_w \pi^\prime,$ then $\beta \circ \pi
\circ \beta \le_w \beta \circ \pi^\prime \circ \beta.$  By combining
Proposition \ref{lengthen} with cases where Conjecture \ref{main
  conjecture} is known to hold many pairs $T,S$ with $S$ dominating
$T$ can be constructed.  For instance, ${\bf abbab}$ dominates ${\bf
  abaab},$ and ${\bf bbbaab}$ dominates ${\bf abbaaa}$ are known cases
of Conjecture \ref{main conjecture}.  Hence, ${\bf abbbaabaaabbab}$
dominates ${\bf aabbaaaaaabaab}.$

If true, Conjecture \ref{main conjecture} would be surprising.  We are
required to find $\phi:D(T) \to D(T \circ \beta) $ such that $I(\pi)
\subseteq I(\phi(\pi)).$ Yet, $I(\pi) \cap I(\pi \circ \beta) =
\emptyset$ for all $\pi.$

\section{From flag $h$-vectors to $h$-vectors}  \label{i to i+1}
  Our goal is to give a combinatorial flag $h$-vector proof of
  (\ref{going up}).  One way to do this would be to construct a
  matching from $(i-1)$-subsets of $[r]$ to $i$-subsets of $r$ such
  that each $(i-1)$-subset is matched to an $i$-subset which dominates
  it.  For example, here is such a matching from 2-subsets of $[6]$ to
  3-subsets of $[6].$
  
$$\begin{array}{ccccc}
\{1,2\} \to \{1,2,3\} & \ &\{2,3\} \to \{2,3,4\} & \ & \{3,5\} \to \{1,3,5\} \\
\{1,3\} \to \{1,3,6\} & \ &\{2,4\} \to \{2,4,6\} & \ & \{3,6\} \to \{2,3,6\} \\
\{1,4\} \to \{1,3,4\} & \ &\{2,5\} \to \{2,3,5\} & \ & \{4,5\} \to \{2,4,5\} \\  
\{1,5\} \to \{1,4,5\} & \ &\{2,6\} \to \{2,5,6\} & \ & \{4,6\} \to \{1,4,6\} \\ 
\{1,6\} \to \{1,5,6\} & \ &\{3,4\} \to \{3,4,6\} & \ & \{5,6\} \to \{4,5,6\} \\ 
\end{array}
$$
\noindent This matching gives a combinatorial flag $h$-vector proof
that $h_2 \le h_3$ for all geometric lattices whose rank is greater
than or equal to $7$.  Table \ref{3 to 4} gives a matching for
$[3]$-sets of $[8]$ to $4$-sets of $[8]$.
  As we will see below, this is not possible for all $r$ and $i \le
  r/2,$ but it can be done for somewhat smaller $i.$  
 
 \begin{lem}
   Let $T$ be an $i$-subset of $[r].$  Then there exists at least
   $\lfloor r- \frac{5}{2} i \rfloor$ supersets of $T$ of cardinality
   $i+1$ which dominate $T.$  
 \end{lem} 
 
\begin{figure} 
\begin{picture}(200,200)(-100,00)
 
  \put (0,190){$\emptyset$}
  \put (3,185){\line(0,-1){15}}
  \put (5,185){\line(4,-1){25}}
  \put (32,173){${\bf b|}$}
  \put (0,162){${\bf a}$}
  \put (3,157){\line(0,-1){15}}
  \put (5, 157){\line(4,-1){25}}
  \put (32,144){${\bf \dot{a}|a}$}
  \put (-2,134){${\bf ab}$}
  \put (5,129){\line(4,-1){25}}
  \put(32,117){${\bf abb|}$}
  \put (3,129){\line(0,-1){15}}
  \put (-5,101){${\bf aba}$}
  \put (5,96){\line(4,-1){25}}
  \put (32,82){${\bf abab|}$}
  \put (3,96){\line(0,-1){15}}
  \put (-8,65){${\bf abaa}$}
  \put (3,60){\line(0,-1){15}}
  \put (5,60){\line(4,-1){25}}
  \put (32,48){${\bf ab\dot{a}\dot{a}|a}$}
  \put (-11,37){${\bf abaab}$}
  \put (3,32){\line(0,-1){15}}
  \put (5,32){\line(4,-1){25}}
  \put (32,20){${\bf abaabb|}$}
  \put (-13,8){${\bf ab\dot{a}aba|}$}
 
\end{picture} 
 \caption{Initial descent patterns}  \label{tree}
\end{figure}
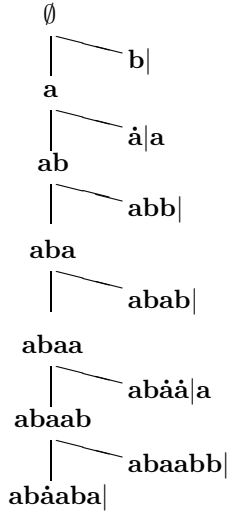

 \begin{proof}
   The proof is by induction on $r,$ the case of $r =1$ being trivial.
   Consider the tree in Figure \ref{tree}.  Each ${\bf ab}$-monomial
   stands for the initial descent pattern of a permutation.  A dot
   above an ascent at position $k$ means that if $T$ has the given
   initial descent pattern, then  $T$ is dominated by $T \cup \{k\}.$
   The bottom of each branch shows how the induction hypothesis should
   be applied.  If the monomial ends with an ascent, then the
   induction hypothesis is applied to everything to the right of the
   bar.  For instance, the entry ${\bf ab\dot{a}\dot{a}|a}$ says to
   apply the induction hypothesis  to the $r-4$ positions to the right
   of ${\bf ab\dot{a}\dot{a}}.$  Indeed, this insures $\lfloor r-4 -
   \frac{5}{2}(i-2) \rfloor + 2 \ge \lfloor r- \frac{5}{2} i \rfloor$
   supersets which dominate $T.$  When the monomial ends with ${\bf
   ...b|}$ the induction hypothesis is applied to the positions to the
   right of the next ascent.  For instance, the monomial ${\bf
   ababbbbaabba}$ is covered by the branch endpoint ${\bf abab|}.$  
   
Let $u=m(T).$  Starting from the top of the tree we look for either an
interior node which matches $u$ exactly or a branch endpoint which
matches an initial segment of $u.$   All of the interior nodes satisfy
the theorem and all of the branch endpoints demonstrate how to use the
induction hypothesis to prove the theorem for $u.$    
   \end{proof}
   
\begin{thm} \label{explained i}
  If $i \le  \frac{2}{7} (r + \frac{5}{2}),$ then there exists a
  matching $\phi$ from $(i-1)$-subsets of $[r]$ to $i$-subsets of
  $[r]$ such that $\phi(T)$ dominates $T$ for each $|T|=i-1.$
\end{thm}

\begin{proof} 
  The condition on $i$ ensures that $r - \frac{5}{2}(i-1) \ge i.$
  Hence, by the above lemma, each $(i-1)$-subset of $[r]$ has at least
  $i$ supersets which dominate it.  Obviously any  $i$-subset of $[r]$
  has at most $i$ subsets of cardinality $i-1$ which it dominates.
  The theorem is now an elementary application of Hall's marriage
  theorem.  
\end{proof}

 Theorem \ref{explained i} is not optimal.  We have already seen that
 there are suitable matchings for 2-sets to 3-sets in $[6]$ and 3-sets
 to 4-sets in $[8].$  However, it is not always possible to obtain
 suitable matchings for all $i \le r/2.$ 

\begin{example}
  Let $T = \{2,5,6,9\} \subseteq [10].$  There are no 5-supersets of
  $T$ in $[10]$ which dominate $T.$  This can be seen by directly
  computing $h_T$ and $h_S$ in $\B_{11}$ where $S$ runs over all
  potential supersets.  In each case $h_T > h_S.$  
\end{example}

\begin{prob}
  Asymptotically, Theorem \ref{explained i} covers a little over 57\%
  of the inequalities in (\ref{going up}).  How much can this be improved?
 \end{prob}
 
 Instead of insisting on a one-to-one matching we can consider
 grouping subsets together.  Of course, this must be done in
 moderation.  Indeed,  $h_{i-1} \le h_i$ is just a reflection of
 grouping all subsets of the same cardinality together.
 
 \begin{example}
   Let $T=\{3\}, S=\{2\}, U= \{2,3\}, V= \{1,3\}.$  As the following
   table shows, there is a bijection which respects the weak Bruhat
   order from the permutations in $S_4$ whose  descent set is $T$ or
   $S$ to those whose descent set is $U$ or $V.$  Hence $h_T + h_S \le
   h_U + h_V$ in rank 4 (or more) geometric lattices.
\end{example}

$$\begin{array}{ccc}
1243 & \to & 2143 \\
1342 & \to & 1432\\
2341 & \to & 2431 \\
1324 & \to & 3142 \\
1423 & \to & 4132 \\
2314 & \to & 3241 \\
2413 & \to & 4231 \\
3412 & \to & 3421
\end{array}$$

\noindent Combined with the previously shown $h_{\{1\}} \le
h_{\{1,2\}},$ the above example provides a combinatorial proof of $h_1
\le h_2$ for geometric lattices of rank 4 or greater. 

As noted earlier, the ${\bf cd}$-index of an oriented matroid  is a
nonnegative linear combination of the ${\bf ab}$-index of the
associated geometric lattice.  
 
\begin{prob}
  Are there groupings of subsets such that the corresponding flag
  $h$-vector inequality translates to a {\bf cd}-inequality for
  oriented matroids?
\end{prob}

\vspace{1 in}

\begin{table}[!hb]
$$\begin{array}{cc|cc|cc}
T & \phi(T) & T & \phi(T) & T & \phi(T) \\ 
\{1,2,3\} & \{1,2,3,4\} & \{1,2,5\} & \{1,2,4,5\} & \{1,2,6\} & \{1,2,4,6\} \\
\{1,2,8\} & \{1,2,3,8\} & \{1,3,4\} & \{1,3,4,8\} & \{1,3,6\} & \{1,3,6,8\} \\
\{1,3,7\} & \{1,3,6,7\} & \{1,3,8\} & \{1,3,5,8\} & \{1,4,5\} & \{1,4,5,8\} \\
\{1,4,6\} & \{1,4,6,8\} & \{1,4,7\} & \{1,4,5,7\} & \{1,4,8\} &
\{1,4,7,8\} \\
\{1,5,6\} & \{1,2,5,6\} & \{1,5,7\} & \{1,2,5,7\} & \{1,5,8\} & \{1,2,5,8\} \\
\{1,6,7\} & \{1,2,6,7\} & \{1,6,8\} & \{1,2,6,8\} & \{1,7,8\} &
\{1,6,7,8\} \\
\{2,3,4\} & \{2,3,4,5\} & \{2,3,5\} & \{2,3,5,6\} & \{2,3,6\} & \{2,3,6,8\} \\
\{2,3,7\} & \{2,3,4,7\} & \{2,3,8\} & \{2,3,4,8\} & \{2,4,5\} & \{2,4,5,6\} \\
\{2,4,6\} & \{2,4,6,7\} & \{2,4,7\} & \{2,4,5,7\} & \{2,4,8\} & \{2,4,7,8\} \\
\{2,5,6\} & \{2,5,6,8\} & \{2,5,7\} & \{2,3,5,7\} & \{2,5,8\} &
\{2,4,5,8\} \\
\{2,6,7\} & \{2,3,6,7\} & \{2,6,8\} & \{2,4,6,8\} & \{2,7,8\} & \{2,6,7,8\} \\
\{3,4,5\} & \{1,3,4,5\} & \{3,4,6\} & \{1,3,4,6\} & \{3,4,7\} & \{1,3,4,7\} \\
\{3,4,8\} & \{3,4,7,8\} & \{3,5,6\} & \{1,3,5,6\} & \{3,5,7\} &
\{1,3,5,7\} \\
\{3,5,8\} & \{2,3,5,8\} & \{3,6,7\} & \{3,4,6,7\} & \{3,6,8\} &
\{3,5,6,8\} \\
\{3,7,8\} & \{2,3,7,8\} & \{4,5,6\} & \{1,4,5,6\} & \{4,5,7\} & \{3,4,5,7\} \\
\{4,5,8\} & \{3,4,5,8\} & \{4,6,7\} & \{1,4,6,7\} & \{4,6,8\} &
\{3,4,6,8\} \\
\{4,7,8\} & \{4,5,7,8\} & \{5,6,7\} & \{4,5,6,7\} & \{5,6,8\} & \{4,5,6,8\} \\
\{5,7,8\} & \{1,5,7,8\} & \{6,7,8\} & \{5,6,7,8\} & & 
\end{array}$$
\caption{An injection $\phi$ from 3-sets of $[8]$ to 4-sets of $[8]$
  such that $\phi(T)$ dominates $T.$} \label{3 to 4}
\end{table}

\begin{table}[!ht]
$$\begin{array}{ccc}
T \to T \circ \beta & & T \to T \circ \beta \\
 \ & r = 3 & \\
\{1\} \to \{1,2\} & \  & \{3\} \to \{2,3\} \\
&r=5 & \\
 \{1,2\} \to \{1,2,3\} & & \{4,5\} \to \{3,4,5\} \\
 \{1,4\} \to \{1,3,4\} & & \{2,5\} \to \{2,3,5\} \\
 &r=7& \\
 \{1,2,3\} \to \{1,2,3,4\} & & \{5,6,7\} \to \{4,5,6,7\} \\
 \{1,2,5\} \to \{1,2,4,5\} & & \{3,6,7\} \to \{3,4,6,7\} \\
 \{1,3,6\} \to \{1,3,4,6\} & & \{2,5,7\} \to \{2,4,5,7\} \\
 \{1,5,6\} \to \{1,4,5,6\} & & \{2,3,7\} \to \{2,3,4,7\}
\end{array}$$
\caption{Several known examples of $T \circ \beta$ dominating $T.$}  \label{known cases}
\end{table}

{\bf Acknowledgment.} As has frequently been the case throughout our
careers we benefitted from helpful conversations with Louis Billera.
The first author was partly supported by an NSF-VIGRE postdoctoral fellowship.
The second author was partly supported by NSF grant DMS-0245623.


\begin{thebibliography}{1}

\bibitem[Ai]{Ai}
M.~Aigner.
\newblock Whitney numbers.
\newblock In N.~L. White, editor, {\em Combinatorial Geometries},
pages 139 -- 160. Cambridge University Press, 1987.

\bibitem[BK]{BK}
M.~Bayer, and A.~Klapper.
\newblock A new index for polytopes.
\newblock {\em Discrete and Comput. Geom.}, 6: 33--47, 1991.

\bibitem[BS]{BS}
M.~Bayer, and B.~Sturmfels.
\newblock Lawrence polytopes.
\newblock {\em Canad. J. Math}, 42: 62--79, 1990.

\bibitem[BER]{BER}
Louis~J.~Billera, Richard Ehrenborg, and Margaret Readdy.
\newblock The {\bf c-2d}-index of oriented matroids.
\newblock {\em J. Comb. Theory Ser. A}, 80:79 -- 105, 1997.

\bibitem[Bi]{Bi}
G.~Birkhoff.
\newblock Abstract linear dependence in lattices.
\newblock {\em Amer. J. Math.}, 57: 800--804, 1935.
 
\bibitem[Bj1]{Bj2}
A.~Bj{\"o}rner.
\newblock Shellable and Cohen-Macaulay partially ordered sets.
\newblock {\em Trans. Amer. Math. Soc.}, 260: 159--183, 1980.

\bibitem[Bj2]{Bj5}
A.~Bj{\"o}rner.
\newblock On the homology of the geometric lattice.
\newblock {\em Algebra Universalis}, 14:107--128, 1982.

\bibitem[Bj3]{Bj3}
A.~Bj{\"o}rner.
\newblock Orderings of Coxeter groups.
\newblock In C.~Greene, editor, {\em Combinatorics and Algebra}, pages 175 -- 195.
\newblock Contemporary Mathematics Series Vol. 34. Amer. Math. Soc., 1985.

\bibitem[Bj4]{Bj}
A.~Bj{\"o}rner.
\newblock The homology and shellability of matroids and geometric
lattices.
\newblock In N.~L. White, editor, {\em Matroid Applications}, pages 226 --
283.
  Cambridge University Press, 1992.

\bibitem[BLSWZ]{Bj4}
A.~Bj{\"o}rner, M.~Las~Vergnas, B.~Sturmfels, N.~White, G.~Ziegler.
\newblock {\em Oriented Matroids}.
\newblock Cambridge University Press, Cambridge 1993.
 
\bibitem[Ch]{Ch}
M.K.~Chari.
\newblock Two decompositions in topological combinatorics with
applications to
  matroid complexes.
\newblock {\em Trans. Amer. Math. Soc.}, 349:3925--3943, 1997.

\bibitem[DW1]{DW1}
T.A.~Dowling and R.M.~Wilson.
\newblock The slimmest geometric lattices.
\newblock {\em Trans. Amer. Math. Soc.}, 196: 203--215, 1974.

\bibitem[DW2]{DW2}
T.A.~Dowling and R.M.~Wilson.
\newblock Whitney number inequalities for geometric lattices.
\newblock {\em Proc. Amer. Math. Soc.}, 47: 504--512, 1975.

\bibitem[Fo]{Fo}
J.~Folkman.
\newblock The homology groups of a lattice.
\newblock {\em J. Math. Mech.}, 15:631--636, 1966.

\bibitem[Ny]{Ny}
K.L.~Nyman.
\newblock Enumeration in geometric lattices and the symmetric group.
\newblock Ph.D. thesis, Cornell University, Ithaca, NY, 2001.

\bibitem[Ox]{Ox}
J.~Oxley
\newblock {\em Matroid theory}.
\newblock Oxford University Press, 1992.

\bibitem[St]{St}
R.P.~Stanley.
\newblock {\em Combinatorics and commutative algebra}.
\newblock Progress in Mathematics, 41. Birkh{\"{a}user} Boston, Inc.,
1996.

\bibitem[St2]{St3}
R.P.~Stanley.
\newblock {\em Enumerative combinatorics}.
\newblock Vol. 1, Cambridge University Press, Cambridge, 1997.

\bibitem[St3]{St2}
R.P.~Stanley.
\newblock Balanced Cohen-Macaulay complexes.
\newblock {\em Trans. Amer. Math. Soc.}, 249: 139--157, 1979.

\bibitem[Za]{Za}
T.~Zaslavsky.
\newblock The M\"{o}bius function and the characteristic polynomial.
\newblock In N.~L. White, editor, {\em Combinatorial Geometries},
pages 114 -- 138. Cambridge University Press, 1987.


\end{thebibliography}
\end{document}